\documentclass{svjour3}
\usepackage[T1]{fontenc}
\usepackage[latin9]{inputenc}
\usepackage{prettyref}
\usepackage{float}
\usepackage{url}
\usepackage{amsmath}
\usepackage{amssymb}
\usepackage{esint}
\usepackage{xargs}[2008/03/08]
\usepackage[unicode=true,pdfusetitle,
 bookmarks=true,bookmarksnumbered=false,bookmarksopen=false,
 breaklinks=false,pdfborder={0 0 1},backref=false,colorlinks=false]
 {hyperref}
\usepackage{breakurl}

\makeatletter

\floatstyle{ruled}
\newfloat{algorithm}{tbp}{loa}
\providecommand{\algorithmname}{Algorithm}
\floatname{algorithm}{\protect\algorithmname}

\newenvironment{svmultproof}{\begin{proof}}{\qed\end{proof}}

\usepackage{comment}

\usepackage{amssymb,yhmath}

\newcommand{\isdef}{\ensuremath{\stackrel{\text{def}}{=}}}

\newcommand{\ZZ}{\ensuremath{\mathbb{Z}}}

\newcommand{\RR}{\ensuremath{\mathbb{R}}}
\newcommand{\CC}{\ensuremath{\mathbb{C}}}
\newcommand{\der}[2]{#1^{(#2)}}

\DeclareMathOperator{\dd}{d}
\DeclareMathOperator{\ee}{e}
\DeclareMathOperator{\diag}{diag}

\usepackage{prettyref}
\newrefformat{prob}{Problem \ref{#1}}
\newrefformat{prop}{Proposition \ref{#1}}
\newrefformat{def}{Definition \ref{#1}}
\newrefformat{sub}{Section \ref{#1}}
\newrefformat{alg}{Algorithm \ref{#1}}
\newrefformat{con}{Conjecture \ref{#1}}
\newrefformat{rem}{Remark \ref{#1}}

\makeatother

\begin{document}

\title{Local and global geometry of Prony systems and Fourier reconstruction
of piecewise-smooth functions\thanks{This research is supported by the Adams Fellowship Program of the Israel Academy of Sciences and Humanities, ISF grant 264/09 and the Minerva Foundation.}}

\titlerunning{Geometry of Prony systems and reconstruction of piecewise-smooth
functions}

\author{D.Batenkov \and Y.Yomdin}

\institute{Department of Mathematics, The Weizmann Institute of Science, Rehovot
76100, Israel.}
\maketitle
\begin{abstract}
Many reconstruction problems in signal processing require solution
of a certain kind of nonlinear systems of algebraic equations, which
we call Prony systems. We study these systems from a general perspective,
addressing questions of global solvability and stable inversion. Of
special interest are the so-called ``near-singular'' situations,
such as a collision of two closely spaced nodes.

We also discuss the problem of reconstructing piecewise-smooth functions
from their Fourier coefficients, which is easily reduced by a well-known
method of K.Eckhoff to solving a particular Prony system. As we show
in the paper, it turns out that a modification of this highly nonlinear
method can reconstruct the jump locations and magnitudes of such functions,
as well as the pointwise values between the jumps, with the maximal
possible accuracy. 
\end{abstract}
\global\long\def\g{\gamma}
\global\long\def\M{\tilde{M}}
\global\long\def\PM{{\cal PM}}
\global\long\def\err{\varepsilon}
\newcommandx\confvec[2][usedefault, addprefix=\global, 1=j, 2=k]{\vec{v_{#1,#2}}}
 \global\long\def\cvand{\ensuremath{V}}
 \global\long\def\o{\omega}

\section{Introduction}

In many applications, it is often required to reconstruct an unknown
signal from a small number of measurements, utilizing some a-priori
knowledge about the signal structure. Such problems arose (and continue
to arise) in recent years under several names in different fields,
such as Finite Rate of Innovation, super-resolution, sub-Nyquist sampling
and Algebraic Signal Reconstruction \cite{bat2008,BatGolYom2011,BatenkovFDE,candes2012towards,donoho1992superresolution,donoho2006stable,ettinger2008lvn,gustafsson2000rpd,sig_ack,vetterli2002ssf}.
One underlying connection between these problems is that almost all
of them require solution of a certain kind of nonlinear systems of
algebraic equations, which we call Prony systems. We therefore consider
the study of this system to be an important topic. In particular,
questions of solvability, uniqueness, including in near-singular situations,
as well as stability of reconstruction in the presence of noise turn
out to be non-trivial and requiring a delicate study of some related
algebraic-geometric structures.

This paper consists of two parts. In the first part, we consider the
general Prony system. First, we present a necessary and sufficient
condition for the system to be solvable. Next, we give simple estimate
of the stability of inversion in a ``regular'' setting. Finally,
we consider inversion in several ``near-singular'' situations, and
in particular the practically important situation of colliding nodes.
We show that a reparametrization in the basis of divided finite differences
turns the problem into a well-posed one in this setting.

In the second part of the paper, we present our recent solution to
a conjecture posed by K.Eckhoff in 1995 \cite{eckhoff1995arf}, which
asks for an algorithm to reconstruct a piecewise-smooth function with
unknown discontinuity locations from its first Fourier coefficients.
While the problem of defeating the Gibbs phenomenon received much
attention in the last decades (see \cite{adcock2010stable,adcock2012stability,batyomAlgFourier,eckhoff1995arf,eckhoff1998high,gelb1999detection,gottlieb1997gpa,hrycak2010pseudospectral,jerriGibbs11,tadmor2007filters,wei2007detection}
and references therein), the question of attaining maximal possible
accuracy of reconstruction remained open. We show how the Algebraic
Reconstruction approach, and in particular an accurate solution of
a certain Prony system, provides the required approximation rate.

\section{The Prony problem}

Prony system appears as we try to solve a very simple ``algebraic
signal reconstruction\textquotedbl{} problem of the following form:
assume that the signal $F(x)$ is known to be a linear combination
of shifted $\delta$-functions:
\begin{equation}
F\left(x\right)=\sum_{j=1}^{d}a_{j}\delta\left(x-x_{j}\right).\label{eq:equation_model_delta}
\end{equation}
We shall use as measurements the polynomial moments:
\begin{equation}
m_{k}=m_{k}\left(F\right)=\int x^{k}F\left(x\right)\dd x.\label{eq:moments}
\end{equation}
After substituting $F$ into the integral defining $m_{k}$ we get
\[
m_{k}(F)=\int x^{k}\sum_{j=1}^{d}a_{j}\delta(x-x_{j})\dd x=\sum_{j=1}^{d}a_{j}x_{j}^{k}.
\]
Considering $a_{j}$ and $x_{j}$ as unknowns, we obtain equations
\begin{equation}
m_{k}\left(F\right)=\sum_{j=1}^{d}a_{j}x_{j}^{k},\; k=0,1,\dots.\label{eq:equation_prony_system}
\end{equation}
This infinite set of equations (or its part, for $k=0,1,\dots,2d-1$),
is called Prony system. It can be traced at least to R. de Prony (1795,
\cite{prony1795essai}) and it is used in a wide variety of theoretical
and applied fields. See \cite{Auton1981} for an extensive bibligoraphy
on the Prony method.

In writing Prony system \eqref{eq:equation_prony_system} we have
assumed that all the nodes $x_{1},\dots,x_{d}$ are pairwise different.
However, as a right-hand side $\mu=(m_{0},\dots,m_{2d-1})$ of \eqref{eq:equation_prony_system}
is provided by the actual measurements of the signal $F$, we cannot
guarantee a priori, that this condition is satisfied for the solution.
Moreover, we shall see below that multiple nodes may naturally appear
in the solution process. In order to incorporate possible collisions
of the nodes, we consider ``confluent Prony systems''.

Assume that the signal $F(x)$ is a linear combination of shifted
$\delta$-functions and their derivatives:
\begin{equation}
F\left(x\right)=\sum_{j=1}^{s}\sum_{\ell=0}^{d_{j}-1}a_{j,\ell}\delta^{\left(\ell\right)}\left(x-x_{j}\right).\label{eq:confluent equation_model_delta}
\end{equation}

\begin{definition}
\label{def:mult-vec-signal}For $F\left(x\right)$ as above, the vector
$D\left(F\right)=(d_{1},\dots,d_{s})$ is \emph{the multiplicity vector}
of $F$, $s=s\left(F\right)$ is its \emph{degree }and $d=\sum_{j=1}^{s}d_{j}$
is its \emph{order}. For avoiding ambiguity in these definitions,
it is always understood that $a_{j,d_{j}-1}\neq0$ for all $j=1,\dots,s$.
\end{definition}
For the moments $m_{k}=m_{k}(F)=\int x^{k}F(x)\dd x$ we now get 
\[
m_{k}=\sum_{j=1}^{s}\sum_{\ell=0}^{d_{j}-1}a_{j,\ell}\frac{{k!}}{{(k-\ell)!}}x_{j}^{k-\ell}.
\]
Considering $x_{i}$ and $a_{j,\ell}$ as unknowns, we obtain a system
of equations
\begin{equation}
\sum_{j=1}^{s}\sum_{\ell=0}^{d_{j}-1}\frac{k!}{\left(k-\ell\right)!}a_{j,\ell}x_{j}^{k-\ell}=m_{k},\quad k=0,1,\dots,2d-1,\label{eq:confluent equation_prony_system}
\end{equation}
which is called a confluent Prony system of order $d$ with the multiplicity
vector $D=\left(d_{1},\dots,d_{s}\right)$. The original Prony system
\eqref{eq:equation_prony_system} is a special case of the confluent
one, with $D$ being the vector $(1,\dots,1)$ of the length $d$.

The system \eqref{eq:confluent equation_prony_system} arises also
in the problem of reconstructing a planar polygon $P$ (or even an
arbitrary semi-analytic \emph{quadrature domain}) from its moments
\[
m_{k}(\chi_{P})=\iint_{\RR^{2}}z^{k}\chi_{P}\dd x\dd y,\; z=x+\imath y,
\]
where $\chi_{P}$ is the characteristic function of the domain $P\subset\RR^{2}$.
This problem is important in many areas of science and engineering
\cite{gustafsson2000rpd}. The above yields the confluent Prony system
\[
m_{k}=\sum_{j=1}^{s}\sum_{i=0}^{d_{j}-1}c_{i,j}k(k-1)\cdots(k-i+1)z_{j}^{k-i},\qquad c_{i,j}\in\CC,\; z_{j}\in\CC\setminus\left\{ 0\right\} .
\]

As we shall see below, if we start with the measurements $\mu(F)=\mu=(m_{0},\dots,m_{2d-1})$,
then a natural setting of the problem of solving the Prony system
is the following:
\begin{problem}[Prony problem of order $d$]
\label{prob:prony}\textit{Given the measurements 
\[
\mu=(m_{0},\dots,m_{2d-1})\in\CC^{2d}
\]
in the right hand side of \eqref{eq:confluent equation_prony_system},
find the multiplicity vector $D=(d_{1},\dots,d_{s})$ of order $r=\sum_{j=1}^{s}d_{j}\leq d$,
and find the unknowns $x_{j}$ and $a_{j,\ell},$ which solve the
corresponding confluent Prony system \eqref{eq:confluent equation_prony_system}
with the multiplicity vector $D$.}
\end{problem}
It is extremely important in practice to have \emph{a stable method
of inversion}. Many research efforts are devoted to this task (see
e.g. \cite{badeau2008performance,batenkov2011accuracy,donoho2006stable,peter2011nonlinear,potts2010parameter,stoica1989music}
and references therein). A basic question here is the following.
\begin{problem}[Noisy Prony problem]
\label{prob:noisy-prony}Given the \emph{noisy} measurements \textit{
\[
\tilde{\mu}=(\tilde{m_{0}},\dots,\tilde{m}_{2d-1})\in\CC^{2d}
\]
and an estimate of the error $\left|\tilde{m}_{k}-m_{k}\right|\leq\err_{k}$,
solve \prettyref{prob:prony} so as to minimize the reconstruction
error.}
\end{problem}

\section{Solving the Prony problem}

\subsection{Prony mapping}

Let us introduce some notations which will be useful in subsequent
treatment.
\begin{definition}
\label{def:config}For each $w=\left(x_{1},\dots,x_{d}\right)\in\CC^{d}$,
let $s=s\left(w\right)$ be the number of distinct coordinates $\tau_{j}$,
$j=1,\dots,s$, and denote $T\left(w\right)=\left(\tau_{1},\dots,\tau_{s}\right)$.
The multiplicity vector is $D=D\left(w\right)=\left(d_{1},\dots,d_{s}\right)$,
where $d_{j}$ is the number of times the value $\tau_{j}$ appears
in $\left\{ x_{1},\dots,x_{d}\right\} .$ The order of the values
in $T\left(w\right)$ is defined by their order of appearance in $w$.\end{definition}
\begin{example}
For $w=\left(3,1,2,1,0,3,2\right)$ we have $s=4$, $T\left(w\right)=\left(3,1,2,0\right)$
and $D\left(w\right)=\left(2,2,2,1\right)$.\end{example}
\begin{remark}
\label{rem:mult-vec-abuse}Note the slight abuse of notations between
\prettyref{def:mult-vec-signal} and \prettyref{def:config}. Note
also that the \emph{order }of \emph{$D\left(w\right)$ }equals to\emph{
$d$ }\textbf{\emph{for all }}$w\in\CC^{d}$\emph{.}\end{remark}
\begin{definition}
For each $w\in\CC^{d}$, let $s=s\left(w\right),\; T\left(w\right)=\left(\tau_{1},\dots,\tau_{s}\right)$
and $D\left(w\right)=\left(d_{1},\dots,d_{s}\right)$ be as in \prettyref{def:config}.
We denote by $V_{w}$ the vector space of dimension $d$ containing
the linear combinations
\begin{equation}
g=\sum_{j=1}^{s}\sum_{\ell=0}^{d_{j}-1}\g_{j,\ell}\delta^{\left(\ell\right)}\left(x-\tau_{j}\right)\label{eq:standard-basis-representation}
\end{equation}
of $\delta$-functions and their derivatives at the points of $T\left(w\right)$.
The ``standard basis'' of $V_{w}$ is given by the distributions
\begin{equation}
\delta_{j,\ell}=\delta^{\left(\ell\right)}\left(x-\tau_{j}\right),\qquad j=1,\dots,s\left(w\right);\;\ell=0,\dots,d_{j}-1.\label{eq:standard-basis-prony}
\end{equation}

\end{definition}
\begin{minipage}[t]{1\columnwidth}%
\end{minipage}
\begin{definition}
The Prony space ${\cal P}_{d}$ is the vector bundle over $\CC^{d}$,
consisting of all the pairs
\[
\left(w,g\right):\quad w\in\CC^{d},\; g\in V_{w}.
\]
The topology on ${\cal P}_{d}$ is induced by the natural embedding
${\cal P}_{d}\subset\CC^{d}\times{\cal D},$ where ${\cal D}$ is
the space of distributions on $\CC$ with its standard topology.
\end{definition}
Finally, we define the Prony mapping ${\cal PM}$ which encodes the
Prony problem.
\begin{definition}
The Prony mapping $\PM:{\cal P}_{d}\to\CC^{2d}$ for $\left(w,g\right)\in{\cal P}_{d}$
is defined as follows:
\[
{\cal \PM}\left(\left(w,g\right)\right)=\left(m_{0},\dots,m_{2d-1}\right)\in\CC^{2d},\qquad m_{k}=m_{k}\left(g\right)=\int x^{k}g\left(x\right)\dd x.
\]

\end{definition}
Therefore, a formal solution of the Prony problem is given by the
inversion of the Prony mapping $\PM$.

Finally, let us recall an important type of matrices which play a
central role in what follows.
\begin{definition}
Let $\left(x_{1},\dots,x_{s}\right)\in\CC^{s}$ and $D=\left(d_{1},\dots,d_{s}\right)$
with $d=\sum_{j=1}^{s}d_{j}$ be given. The\emph{ $d\times d$ confluent
Vandermonde} matrix is
\begin{equation}
\cvand=\cvand\left(x_{1},d_{1},\dots,x_{s},d_{s}\right)=\left[\begin{array}{cccc}
\confvec[1][0] & \confvec[2][0] & \dotsc & \confvec[s][0]\\
\confvec[1][1] & \confvec[2][1] & \dotsc & \confvec[s][1]\\
 &  & \dotsc\\
\confvec[1][d-1] & \confvec[2][d-1] & \dotsc & \confvec[s][d-1]
\end{array}\right]\label{eq:confluent-vandermonde-def}
\end{equation}
where the symbol $\confvec$ denotes the following $1\times d_{j}$
row vector
\[
\confvec\isdef\left[\begin{array}{cccc}
x_{j}^{k}, & kx_{j}^{k-1}, & \dots & ,k\left(k-1\right)\cdots\left(k-d_{j}\right)x_{j}^{k-d_{j}+1}\end{array}\right].
\]

\end{definition}
The matrix $\cvand$ defines the linear part of the confluent Prony
system \eqref{eq:confluent equation_prony_system}, namely, 
\begin{equation}
\cvand\left(x_{1},d_{1},\dots,x_{s},d_{s}\right)\begin{bmatrix}a_{1,0}\\
\vdots\\
a_{1,d_{1}-1}\\
\vdots\\
\\
a_{s,d_{s}-1}
\end{bmatrix}=\begin{bmatrix}m_{0}\\
m_{1}\\
\vdots\\
\\
\\
m_{d-1}
\end{bmatrix}.\label{eq:confluent-prony-linear-part}
\end{equation}

\subsection{Padé problem and the solvability set}

It can be shown that the solution to \prettyref{prob:prony} is equivalent
to solving the well-known Padé approximation problem. While this connection
is extremely important and insightful, we do not provide the details
here for the sake of brevity. Let us only mention the following result.
\begin{proposition}
\label{prop:pade-correspondence}The tuple 
\[
\left\{ s,\; D=(d_{1},\dots,d_{s}),\; r=\sum_{j=1}^{s}d_{j}\leq d,\; X=\left\{ x_{j}\right\} _{j=1}^{s},\; A=\left\{ a_{j,\ell}\right\} _{j=1,\dots,s;\;\ell=0,\dots,d_{j}-1}\right\} 
\]
is a solution to \prettyref{prob:prony} with right-hand side 
\[
\mu=(m_{0},\dots,m_{2d-1})\in\CC^{2d}
\]
if and only if $\left(m_{0},\dots,m_{2d-1}\right)$ are the first
$2d$ Taylor coefficients at $z=\infty$ of the rational function
\[
R_{D,X,A}\left(z\right)=\sum_{j=1}^{s}\sum_{\ell=1}^{d_{j}}\left(-1\right)^{\ell-1}\left(\ell-1\right)!\frac{a_{j,\ell}}{\left(z-x_{j}\right)^{\ell}}=\sum_{k=0}^{2d-1}\frac{m_{k}}{z^{k+1}}+O\left(z^{-2d-1}\right).
\]
The function $R_{D,X,A}\left(z\right)$ is the Stieltjes transform
of the corresponding signal $F\left(x\right)=\sum_{j=1}^{s}\sum_{\ell=0}^{d_{j}-1}a_{j,\ell}\delta^{\left(\ell\right)}\left(x-x_{j}\right)$,
i.e.
\[
R_{D,X,A}\left(z\right)=\int_{-\infty}^{\infty}\frac{F\left(x\right)\dd x}{z-x}.
\]

\end{proposition}
Using this correspondence, it is not difficult to prove the following
result (see \cite{byPronySing12}).
\begin{theorem}
Let the right-hand side $\left(m_{0},\dots,m_{2d-1}\right)$ of \prettyref{prob:prony}
be given. Let $\M_{d}$ denote the $d\times\left(d+1\right)$ Hankel
matrix
\[
\M_{d}=\begin{bmatrix}m_{0} & m_{1} & m_{2} & \dots & m_{d}\\
m_{1} & m_{2} & m_{3} & \dots & m_{d+1}\\
\adots & \adots & \adots & \adots & \adots\\
m_{d-1} & m_{d} & m_{d+1} & \dots & m_{2d-1}
\end{bmatrix}.
\]
For each $e\leqslant d$, denote by $\M_{e}$ the $e\times\left(e+1\right)$
submatrix of $\M_{d}$ formed by the first $e$ rows and $e+1$ columns,
and let $M_{e}$ denote the corresponding square matrix.

Let $r\leqslant d$ be the rank of $\M_{d}$. Then \prettyref{prob:prony}
is solvable if and only if the upper left minor $\left|M_{r}\right|$
of $\M_{d}$ is non-zero. The solution, if it exists, is unique, up
to a permutation of the nodes $\left\{ x_{j}\right\} $. The multiplicity
vector $D=\left(d_{1},\dots,d_{s}\right)$, $\sum_{j=1}^{s}d_{j}=r$,
of the resulting confluent Prony system of order $r$ is the multiplicity
vector of the poles of the rational function $R_{D,X,A}\left(z\right)$,
solving the Padé problem in \prettyref{prop:pade-correspondence}.
\end{theorem}
As a corollary we get a complete description of the right-hand side
data $\mu\in\CC^{2d}$ for which the Prony problem is solvable (unsolvable).
Define for $r=1,\dots,d$ sets $\Sigma_{r}\subset\CC^{2d}$ (respectively,
$\Sigma'_{r}\subset\CC^{2d}$) consisting of $\mu\in\CC^{2d}$ for
which the rank of $\tilde{M}_{d}=r$ and $|M_{r}|\ne0$ \ (respectively,
$|M_{r}|=0)$. The set $\Sigma_{r}$ is a difference $\Sigma_{r}=\Sigma_{r}^{1}\setminus\Sigma_{r}^{2}$
of two algebraic sets: $\Sigma_{r}^{1}$ is defined by vanishing of
all the $s\times s$ minors of $\tilde{M}_{d},\ r<s\leq d,$ while
$\Sigma_{r}^{2}$ is defined by vanishing of $|M_{r}|.$ In turn,
$\Sigma'_{r}=\Sigma_{r}^{'1}\setminus\Sigma_{r}^{'2},$ with $\Sigma_{r}^{'1}=\Sigma_{r}^{1}\cap\Sigma_{r}^{2}$
and $\Sigma_{r}^{'2}$ defined by vanishing of all the $r\times r$
minors of $\tilde{M}_{d}.$ The union $\Sigma_{r}\cup\Sigma'_{r}$
consists of all $\mu$ for which the rank of $\tilde{M}_{d}=r,$ which
is $\Sigma_{r}^{1}\setminus\Sigma_{r}^{'2}.$
\begin{corollary}
\label{cor:non.solv.set}The set $\Sigma$ (respectively, $\Sigma'$)
of $\mu\in\CC^{2d}$ for which the Prony problem is solvable (respectively,
unsolvable) is the union $\Sigma=\cup_{r=1}^{d}\Sigma_{r}$ (respectively,
$\Sigma'=\cup_{r=1}^{d}\Sigma'_{r}$). In particular, $\Sigma'\subset\{\mu\in\CC^{2d},\det M_{d}=0\}.$
\end{corollary}
So for a generic right hand side $\mu$ we have $|M_{d}|\ne0$, and
the Prony problem is solvable. On the algebraic hypersurface of $\mu$
for which $|M_{d}|=0,$ the Prony problem is solvable if $M_{d-1}\ne0$,
etc.

\subsection{\label{sub:stable-inv-away}Stable inversion away from singularities}

\global\long\def\rpm{{\cal \PM}_{D_{0}}^{*}}

Consider \prettyref{prob:noisy-prony} at some interior point $\mu_{0}\in\Sigma$.
By definition, $\mu_{0}\in\Sigma_{r_{0}}$ for some $r_{0}\leq d$.
Let $\left(w_{0},g_{0}\right)={\cal PM}^{-1}\left(\mu_{0}\right)$.
Assume for a moment that the multiplicity vector $D_{0}=D\left(g_{0}\right)=\left(d_{1},\dots d_{s_{0}}\right)$,
$\sum_{j=1}^{s_{0}}d_{j}=r_{0}$, has a non-trivial collision pattern,
i.e. $d_{j}>1$ for at least one $j=1,\dots,s_{0}$. It means, in
turn, that the function $R_{D_{0},X,A}\left(z\right)$ has a pole
of multiplicity $d_{j}$. Evidently, there exists an arbitrarily small
perturbation $\tilde{\mu}$ of $\mu_{0}$ for which this multiple
pole becomes a cluster of single poles, thereby changing the multiplicity
vector to some $D'\neq D_{0}$. While we address this problem in \prettyref{sec:singularities}
via the bases of divided differences, in this section we consider
a ``restricted'' Prony problem. 
\begin{definition}
Let $\PM\left(w_{0},g_{0}\right)=\mu_{0}\in\Sigma_{r_{0}}$ with $D\left(g_{0}\right)=D_{0}$
and $s\left(g_{0}\right)=s_{0}$. Let ${\cal P}_{D_{0}}$ denote the
following subbundle of ${\cal P}_{d}$ of dimension $s_{0}+r_{0}$:
\[
{\cal P}_{D_{0}}=\left\{ \left(w,g\right)\in{\cal P}_{d}:\quad D\left(g\right)=D_{0}\right\} .
\]

The restricted Prony mapping $\rpm:{\cal P}_{D_{0}}\to\CC^{s_{0}+r_{0}}$
is the composition 
\[
\rpm=\pi\circ\PM\restriction_{{\cal P}_{D_{0}}},
\]
where $\pi:\CC^{2d}\to\CC^{s_{0}+r_{0}}$ is the projection map on
the first $s_{0}+r_{0}$ coordinates.
\end{definition}
Inverting this $\rpm$ represents the solution of the confluent Prony
system \eqref{eq:confluent equation_prony_system} with fixed structure
$D_{0}$ from the first $k=0,1,\dots,s_{0}+r_{0}-1$ measurements.
\begin{theorem}[\cite{batenkov2011accuracy}]
Let $\mu_{0}^{*}=\rpm\left(\left(w_{0},g_{0}\right)\right)\in\CC^{s_{0}+r_{0}}$
with the unperturbed solution $g_{0}=\sum_{j=1}^{s_{0}}\sum_{\ell=0}^{d_{j}-1}a_{j,\ell}\delta^{\left(\ell\right)}\left(x-\tau_{j}\right)$.
In a small neighborhood of $\left(w_{0},g_{0}\right)\in{\cal P}_{D_{0}}$,
the map $\rpm$ is invertible. Consequently, for small enough $\err$,
the restricted Prony problem with input data $\tilde{\mu}^{*}\in\CC^{r_{0}+s_{0}}$
satisfying $\|\tilde{\mu}^{*}-\mu_{0}^{*}\|\leq\err$ has a unique
solution. The error in this solution satisfies
\begin{eqnarray*}
\left|\Delta a_{j,\ell}\right| & \leq & \frac{2}{\ell!}\left(\frac{2}{\delta}\right)^{s_{0}+r_{0}}\left(\frac{1}{2}+\frac{s_{0}+r_{0}}{\delta}\right)^{d_{j}-\ell}\left(1+\frac{\left|a_{j,\ell-1}\right|}{\left|a_{j,d_{j}-1}\right|}\right)\err,\\
\left|\Delta\tau_{j}\right| & \leq & \frac{2}{d_{j}!}\left(\frac{2}{\delta}\right)^{s_{0}+r_{0}}\frac{1}{\left|a_{j,d_{j}-1}\right|}\err,
\end{eqnarray*}
where $\delta\isdef\min_{i\neq j}\left|\tau_{i}-\tau_{j}\right|$
(for consistency we take $a_{j,-1}=0$ in the above formula).\end{theorem}
\begin{svmultproof}
[outline]The Jacobian of $\rpm$ can be easily computed, and it
turns out to be equal to the product
\[
{\cal J}_{\rpm}=V\left(\tau_{1},d_{1}+1,\dots,\tau_{s_{0}},d_{s_{0}}+1\right)\diag\left\{ E_{j}\right\} 
\]
where $V$ is the \emph{confluent Vandermonde matrix \eqref{eq:confluent-vandermonde-def}}
on the nodes $\left(\tau_{1},\dots,\tau_{s_{0}}\right)$ and multiplicity
vector
\[
\tilde{D}_{0}=\left(d_{1}+1,\dots,d_{s_{0}}+1\right),
\]
while $E$ is the $\left(d_{j}+1\right)\times\left(d_{j}+1\right)$
block
\[
E_{j}=\begin{bmatrix}1 & 0 & 0 & \cdots & 0\\
0 & 1 & 0 & \cdots & a_{j,0}\\
\vdots & \vdots & \vdots & \ \ddots & \vdots\\
0 & 0 & 0 & \cdots & a_{j,d_{j}-1}
\end{bmatrix}.
\]
Since $\mu_{0}\in\Sigma_{r}$, the highest order coefficients $a_{j,d_{j}-1}$
are nonzero. Furthermore, since all the $\tau_{j}$ are distinct,
the matrix $V$ is nonsingular. Local invertability follows. To estimate
the norm of the inverse, use bounds from \cite{batPronyDec}.
\end{svmultproof}
Let us stress that we are not aware of any general method of inverting
$\rpm$, i.e. solving the restricted confluent Prony problem with
the smallest possible number of measurements. As we shall see below
in \prettyref{sec:eckhoff}, such a method exists for a very special
case of a single point, i.e. $s=1$.

\section{\label{sec:singularities}Prony inversion near singularities}

\subsection{Collision singularities and finite differences}

Collision singularities occur in Prony systems as some of the nodes
$x_{i}$ in the signal $F(x)=\sum_{i=1}^{d}a_{i}\delta(x-x_{i})$
approach one another. This happens for $\mu$ near the discriminant
stratum $\Delta\subset\CC^{2d}$ consisting of those $(m_{0},\dots,m_{2d-1})$
for which some of the coordinates $\left\{ x_{j}\right\} $ in the
solution collide, i.e. the function $R_{D,X,A}\left(z\right)$ has
multiple poles (or, nontrivial multiplicity vector $D$). As we shall
see below, typically, as $\mu$ approaches $\mu_{0}\in\Delta$, i.e.
some of the nodes $x_{i}$ collide, the corresponding coefficients
$a_{i}$ tend to infinity. Notice, that all the moments $m_{k}=m_{k}(F)$
remain bounded. This behavior creates serious difficulties in solving
``near-colliding'' Prony systems, both in theoretical and practical
settings. Especially demanding problems arise in the presence of noise.
The problem of improvement of resolution in reconstruction of colliding
nodes from noisy measurements appears in a wide range of applications.
It is usually called a ``super-resolution problem'' and a lot of
recent publications are devoted to its investigation in various mathematical
and applied settings. See \cite{candes2012towards} and references
therein for a very partial sample.

Here we continue our study of collision singularities in Prony systems,
started in \cite{yom2009Singularities}. The full details will be
published in \cite{byPronySing12}. Our approach uses bases of finite
differences in the Prony space ${\cal P}_{d}$ in order to ``resolve''
the linear part of collision singularities. In these bases the coefficients
do not blow up any more, as some of the nodes collide.

Let $\mu_{0}\in\Sigma_{d}$. Consider the noisy Prony problem in a
neighborhood of the exact solution $\left(w_{0},g_{0}\right)=\PM^{-1}\left(\mu_{0}\right)$.
As explained in \prettyref{sub:stable-inv-away}, if $D\left(w_{0}\right)$
is non-trivial, then there will always be a multiplicity-destroying
perturbation, no matter how small a neighborhood. Assume that the
node vector $w=w\left(\tilde{\mu}\right)$ is determined, and consider
the linear system \eqref{eq:confluent-prony-linear-part} for recovering
the coefficients $\vec{a}=\left\{ a_{j,\ell}\right\} $ in the standard
basis \eqref{eq:standard-basis-prony} of $V_{w\left(\tilde{\mu}\right)}$.
As $\tilde{\mu}\to\mu_{o}$, the matrix of this linear system will
be $\cvand\left(w\left(\tilde{\mu}\right)\right)$ with collision
pattern $D\left(w\left(\tilde{\mu}\right)\right)\neq D_{0}$, and
therefore its determinant will generically approach zero. This will
make the determination of $\left\{ a_{j,\ell}\right\} $ ill-conditioned,
and in fact some of its components will go to infinity. At the limit,
however, the confluent problem is completely well-posed since the
matrix $\cvand\left(w_{0},D_{0}\right)$ is non-singular. The challenge
is, therefore, to make the solution to depend continuously on $\tilde{\mu}$
by a suitable \emph{change of basis} for $V_{w}$. So, instead of
the ill-conditioned system
\[
\tilde{\mu}=V\left(w\left(\tilde{\mu}\right),D\left(w\left(\tilde{\mu}\right)\right)\right)\vec{a}\left(\tilde{\mu}\right)
\]
we would like to have 
\begin{equation}
\tilde{\mu}=\tilde{V}\left(\tilde{\mu}\right)\vec{b}\left(\tilde{\mu}\right),\label{eq:well-posed-linear-recovery}
\end{equation}
where the matrix $\tilde{\cvand}$ is nonsingular and depends continuously
on $\tilde{\mu}$ in the neighborhood of $\mu_{0}$.

First, we extend the well-known definition of divided finite differences
to colliding configurations.
\begin{definition}
\label{def:dd-def}Let $w=\left(x_{1},\dots,x_{d}\right)$ be given.
For each $m=1,2,\dots,d$, denote $w_{m}=\left(x_{1},\dots,x_{m}\right)$.
According to \prettyref{def:config}, let $s_{m}=s\left(w_{m}\right)$,
$T\left(w_{m}\right)=\left(\tau_{1,m},\dots,\tau_{s_{m},m}\right)$
and $D\left(w_{m}\right)=\left(d_{1,m},\dots,d_{s_{m},m}\right)$.
Consider the decomposition of the rational function
\[
R_{w,m}\left(z\right)=\prod_{j=1}^{s_{m}}\frac{1}{\left(z-\tau_{j,m}\right)^{d_{j,m}}}
\]
into the sum of elementary fractions
\begin{equation}
R_{w,m}\left(z\right)=\sum_{j=1}^{s_{m}}\sum_{\ell=1}^{d_{j,m}}\frac{w_{j,l}^{\left(m\right)}}{\left(z-\tau_{j,m}\right)^{\ell}}.\label{eq:dd-decomp-def}
\end{equation}

The $m$-th\emph{ finite difference }$\Delta_{m}\left(w\right)$ is
the following element of $V_{w}$:
\[
\Delta_{m}\left(w\right)=\sum_{j=1}^{s_{m}}\sum_{\ell=1}^{d_{j,m}}\frac{w_{j,\ell}^{\left(m\right)}}{\left(\ell-1\right)!}\delta^{\left(\ell-1\right)}\left(x-\tau_{j,m}\right),
\]
with the coefficients $\left\{ w_{j,\ell}^{\left(m\right)}\right\} $
defined by \eqref{eq:dd-decomp-def}.
\end{definition}
We prove the following results in \cite{byPronySing12}.
\begin{proposition}
\label{prop:continuous-section}The finite difference $\Delta_{m}\left(w\right)$
is a continuous section of the bundle ${\cal P}_{d}$. For $w\in\CC^{d}$
with pairwise distinct coordinates, $\Delta_{m}\left(w\right)$ is
the usual divided finite difference on the elements of $w_{m}$.\end{proposition}
\begin{theorem}
For each $w\in\CC^{d}$, the collection
\[
{\cal B}\left(w\right)=\left\{ \Delta_{m}\left(w\right)\right\} _{m=1}^{d}\subset V_{w}
\]
forms a basis for $V_{w}$.\end{theorem}
\begin{remark}
Another possible way to construct a good basis $\tilde{{\cal B}}\left(w\right)$
is to build the matrix $\tilde{\cvand}$ in \eqref{eq:well-posed-linear-recovery}
directly by imitating the confluence process of the Vandermonde matrices
(as done in \cite{gautschi1962iva}), multiplying $\cvand$ by an
appropriate ``divided difference matrix'' $F\left(\tilde{\mu}\right)$
(more precisely, a chain of such matrices derived from the confluence
pattern). That is, $\tilde{\cvand}=\cvand F$ is the new matrix for
the recovery of the linear part in \eqref{eq:well-posed-linear-recovery},
while the new coefficient vector is $\vec{b}=F^{-1}\vec{a}$. Also
in this case $\tilde{V}\to\cvand\left(w_{0},D_{0}\right)$ as $\tilde{\mu}\to\mu_{0}$.
The matrix $F$ thus defines the corresponding change of basis from
$\left\{ \delta_{j,\ell}\left(w\right)\right\} $ as in \eqref{eq:standard-basis-prony}
to $\tilde{{\cal B}}\left(w\right)$.
\end{remark}
Let us now consider the Prony problem in the basis ${\cal B}\left(w\right)$
in some neighborhood of $\mu_{0}\in\Sigma_{d}$ (thus, the order of
the exact solution $\left(w_{0},g_{0}\right)=\PM^{-1}\left(\mu_{0}\right)$
is $d$). Writing the unknown $g\in V_{w}$ in this basis we have
\[
g=\sum_{m=1}^{d}\beta_{m}\Delta_{m}\left(w\right).
\]

\begin{theorem}[\cite{byPronySing12}]
\label{thm:f.d.prony.1}For $\tilde{\mu}$ in a sufficiently small
neighborhood of $\mu_{0}$, the solution
\[
\PM^{-1}\left(\tilde{\mu}\right)=\left(w\left(\tilde{\mu}\right),\;\left\{ \beta_{m}\left(\tilde{\mu}\right)\right\} \right),
\]
expressed in the basis ${\cal B}\left(w\right)$ of finite differences,
is provided by continuous algebraic functions of $\tilde{\mu}.$\end{theorem}
\begin{svmultproof}
[outline]For each $w$ in a neighborhood of $w_{0}$, we obtain
the system of equations
\begin{equation}
\sum_{m=1}^{d}\beta_{m}\int x^{k}\Delta_{m}\left(w\right)=\tilde{m}_{k},\; k=0,1,\dots,d-1.\label{eq:dd-system-final}
\end{equation}
In the process of solution, the points $\left\{ x_{1},\dots,x_{d}\right\} $
are found as the roots of the polynomial $Q\left(z\right)$ which
appears in the denominator of $R_{X,D,A}\left(z\right)=\frac{P\left(z\right)}{Q\left(z\right)}$.
The coefficient vector $\vec{q}$ of $Q\left(z\right)$ is provided
by solving a non-degenerate linear system
\[
M_{d}\vec{q}=\begin{bmatrix}m_{d}\\
m_{d+1}\\
\vdots\\
m_{2d-1}
\end{bmatrix}.
\]
Therefore, $w=w\left(\tilde{\mu}\right)$ is given by continuous algebraic
functions of $\tilde{\mu}$. By \prettyref{prop:continuous-section},
the functions
\[
\nu_{k,m}\left(w\right)=\int x^{k}\Delta_{m}\left(w\right)\dd x
\]
 are continuous in $w$. At $w=w_{0}$, the system \eqref{eq:dd-system-final}
is non-degenerate by assumption, therefore it stays non-degenerate
in a small neighborhood of $w_{0}$. Thus, the coefficients $\beta_{m}\left(w\left(\tilde{\mu}\right)\right)$
are also continous algebraic functions of $\tilde{\mu}$.
\end{svmultproof}

\subsection{Prony Inversion near $\Sigma'$ and Lower Rank Strata}

The behavior of the inversion of the Prony mapping near the unsolvability
stratum $\Sigma'$ and near the strata where the rank of $\tilde{M}_{d}$
drops, turns out to be pretty complicated. In particular, in the first
case at least one of the nodes tends to infinity. In the second case,
depending on the way the right-hand side $\mu$ approaches the lower
rank strata, the nodes may remain bounded, or some of them may tend
to infinity. In this section we provide one initial result in this
direction, as well as some examples. A comprehensive description of
the inversion of the Prony mapping near $\Sigma'$ and near the lower
rank strata is important both in theoretical study and in applications
of Prony-like systems, and we plan to provide further results in this
direction separately.
\begin{theorem}
\label{thm:near.unsolv}As the right-hand side $\mu\in\CC^{2d}\setminus\Sigma'$
approaches a finite point $\mu_{0}\in\Sigma',$ at least one of the
nodes $x_{1},\dots,x_{d}$ in the solution tends to infinity.\end{theorem}
\begin{svmultproof}
By assumptions, the components $m_{0},\dots,m_{2d-1}$ of the right-hand
side $\mu=(m_{0},\dots,m_{2d-1})\in\CC^{2d}$ remain bounded as $\mu\rightarrow\mu_{0}$.
By \prettyref{thm:f.d.prony.1}, the finite differences coordinates
of the solution $\PM^{-1}(\mu)$ remain bounded as well. Now, if all
the nodes are also bounded, by compactness we conclude that $\PM^{-1}(\mu)\rightarrow\o\in{\cal P}_{d}.$
By continuity in the distribution space (\prettyref{prop:continuous-section})
we have $\PM(\o)=\mu_{0}$. Hence the Prony problem with the right-hand
side $\mu_{0}$ has a solution $\o\in{\cal P}_{d},$ in contradiction
with the assumption that $\mu_{0}\in\Sigma'$.
\end{svmultproof}
As it was shown above, for a given $\mu\in\Sigma$ (say, with pairwise
different nodes) the rank of the matrix $\tilde{M}_{d}$ is equal
to the number of the nodes in the solution for which the corresponding
$\delta$-function enters with a non-zero coefficients. So $\mu$
approaches a certain $\mu_{0}$ belonging to a stratum of a lower
rank of $\tilde{M}_{d}$ if and only if some of the coefficients $a_{j}$
in the solution tend to zero. We do not analyze all the possible scenarios
of such a degeneration, noticing just that if $\mu_{0}\in\Sigma',$
i.e., the Prony problem is unsolvable for $\mu_{0}$, then \prettyref{thm:near.unsolv}
remains true, with essentially the same proof. So at least one of
the nodes, say, $x_{j},$ escapes to infinity. Moreover, one can show
that $a_{j}x_{j}^{2d-1}$ cannot tend to zero - otherwise the remaining
linear combination of $\delta$-functions would provide a solution
for $\mu_{0}$.

If $\mu_{0}\in\Sigma,$ i.e., the Prony problem is solvable for $\mu_{0},$
all the nodes may remain bounded, or some $x_{j}$ may escape to infinity,
but in such a way that $a_{j}x_{j}^{2d-1}$ tends to zero.

\section{\label{sec:eckhoff}Resolution of Eckhoff's problem}

\global\long\def\np{\ensuremath{K}}
\global\long\def\jp{\ensuremath{\xi}}
\global\long\def\jc{\ensuremath{a}}
\newcommandx\fc[1][usedefault, addprefix=\global, 1=k]{\ensuremath{c_{#1}}}
\global\long\def\sc{\ensuremath{M}}
\global\long\def\fun{\ensuremath{f}}
\newcommandx\er[1][usedefault, addprefix=\global, 1=k]{\delta_{#1}}
\global\long\def\err{\varepsilon}
\newcommandx\meas[1][usedefault, addprefix=\global, 1=k]{m_{#1}}
\global\long\def\nmeas{S}
\newcommandx\apprmeas[1][usedefault, addprefix=\global, 1=k]{\widetilde{m}_{#1}}
\newcommandx\frsum[2][usedefault, addprefix=\global, 1=\fun, 2=\sc]{\mathfrak{F}_{#2}\left(#1\right)}
\global\long\def\smooth{\ensuremath{\Psi}}
\global\long\def\sing{\ensuremath{\Phi}}
\global\long\def\nn#1{\widetilde{#1}}
\global\long\def\ord{\ensuremath{d}}
\global\long\def\jcc{\alpha}
\global\long\def\scc{N}
\global\long\def\w{\ensuremath{\omega}}

Consider the problem of reconstructing an integrable function $\fun:\left[-\pi,\pi\right]\to\RR$
from a finite number of its Fourier coefficients
\[
\fc(\fun)\isdef\frac{1}{2\pi}\int_{-\pi}^{\pi}\fun(t)\ee^{-\imath kt}\dd t,\qquad k=0,1,\dots\sc.
\]

It is well-known that for periodic smooth functions, the truncated
Fourier series
\[
\frsum\isdef\sum_{|k|=0}^{\sc}\fc(\fun)\ee^{\imath kx}
\]
converges to $\fun$ very fast, subsequently making Fourier analysis
very attractive in a vast number of applications. By the classical
Jackson's and Lebesgue's theorems \cite{Natanson1949}, if $\fun$
has $d$ continuous derivatives in $\left[-\pi,\pi\right]$ (including
at the endpoints) and $\der{\fun}{\ord}(x)\in Lip\left(R\right)$,
then

\begin{equation}
\max_{-\pi\leq x\leq\pi}\left|\fun(x)-\frsum(x)\right|\leq C\left(R,\ord\right)\sc^{-\ord-1}\ln\sc.\label{eq:best-approximation-smooth}
\end{equation}

Yet many realistic phenomena exhibit discontinuities, in which case
the unknown function $\fun$ is only piecewise-smooth. As a result,
the trigonometric polynomial $\frsum$ no longer provides a good approximation
to $\fun$ due to the slow convergence of the Fourier series (one
of the manifestations of this fact is commonly known as the ``Gibbs
phenomenon''). It has very serious implications, for example when
using spectral methods to calculate solutions of PDEs with shocks.
Therefore an important question arises: \emph{``Can such piecewise-smooth
functions be reconstructed from their Fourier measurements, with accuracy
which is comparable to the 'classical' one \eqref{eq:best-approximation-smooth}''?}

It has long been known that the key problem for Fourier series acceleration
is the detection of the shock locations. Applying elementary considerations
we have the following fact \cite{batFullFourier}.
\begin{proposition}
\label{prop:maximal-accuracy}Let $\fun$ be piecewise $\ord$-smooth.
Then no deterministic algorithm can restore the locations of the discontinuities
from the first $\sc$ Fourier coefficients with accuracy which is
asymptotically higher than $\sc^{-\ord-2}$.
\end{proposition}
Let us first briefly describe what has become known as the Eckhoff's
method for this problem \cite{eckhoff1993,eckhoff1995arf,eckhoff1998high}.

Assume that $\fun$ has $\np>0$ jump discontinuities $\left\{ \jp_{j}\right\} _{j=1}^{\np}$
(they can be located also at $\pm\pi$, but not necessarily so). Furthermore,
we assume that $\fun\in C^{\ord}$ in every segment $\left(\jp_{j-1},\jp_{j}\right)$,
and we denote the associated jump magnitudes at $\jp_{j}$ by
\[
\jc_{\ell,j}\isdef\der{\fun}{\ell}(\jp_{j}^{+})-\der{\fun}{\ell}(\jp_{j}^{-}).
\]
We write the piecewise smooth $\fun$ as the sum $\fun=\smooth+\sing$,
where $\smooth(x)$ is smooth and periodic and $\sing(x)$ is a piecewise
polynomial of degree $\ord$, uniquely determined by $\left\{ \jp_{j}\right\} ,\left\{ \jc_{\ell,j}\right\} $
such that it ``absorbs'' all the discontinuities of $\fun$ and its
first $\ord$ derivatives. This idea is very old and goes back at
least to A.N.Krylov (\cite{kantokryl62}). Eckhoff derives the following
explicit representation for $\sing(x)$: 
\begin{equation}
\begin{split}\sing(x) & =\sum_{j=1}^{\np}\sum_{\ell=0}^{\ord}\jc_{\ell,j}V_{\ell}(x;\jp_{j})\\
V_{n}\left(x;\jp_{j}\right) & =-\frac{\left(2\pi\right)^{n}}{\left(n+1\right)!}B_{n+1}\left(\frac{x-\jp_{j}}{2\pi}\right)\qquad\jp_{j}\leq x\leq\jp_{j}+2\pi
\end{split}
\label{eq:sing-part-explicit-bernoulli}
\end{equation}
where $V_{n}\left(x;\jp_{j}\right)$ is understood to be periodically
extended to $\RR$ and $B_{n}(x)$ is the $n$-th Bernoulli polynomial.
Elementary integration by parts gives the following formula.
\begin{proposition}
Let $\sing(x)$ be given by \eqref{eq:sing-part-explicit-bernoulli}.
For definiteness, let us assume that $\fc[0](\sing)=\int_{-\pi}^{\pi}\sing(x)\dd x=0$.
Then
\begin{equation}
\fc(\sing)=\frac{1}{2\pi}\sum_{j=1}^{\np}\ee^{-\imath k\jp_{j}}\sum_{\ell=0}^{\ord}(\imath k)^{-\ell-1}\jc_{\ell,j},\qquad k=1,2,\dots.\label{eq:singular-fourier-explicit}
\end{equation}

\end{proposition}
Eckhoff observed that if $\smooth$ is sufficiently smooth, then the
contribution of $\fc(\smooth)$ to $\fc(\fun)$ is negligible \textbf{for
large} $k$, and therefore one can hope to reconstruct the unknown
parameters $\left\{ \jp_{j},\jc_{\ell,j}\right\} $ from the perturbed
equations \eqref{eq:singular-fourier-explicit}, where the left-hand
side reads $\fc\left(\fun\right)\sim\fc\left(\sing\right)$ and $k\gg1$.
His proposed method was to construct from the known values
\[
\left\{ \fc\left(\fun\right)\right\} \qquad k=\sc-\left(\ord+1\right)\np+1,\sc-\left(\ord+1\right)\np+2,\dots,\sc
\]
an algebraic equation satisfied by the jump points $\left\{ \jp_{1},\dots,\jp_{\np}\right\} $,
and solve this equation numerically. Based on some explicit computations
for $\ord=1,2;\;\np=1$ and large number of numerical experiments,
he conjectured that his method would reconstruct the jump locations
with accuracy $\sc^{-\ord-1}$. 

We consider the following generalized formulation (without referring
to a specific method).
\begin{conjecture}[Eckhoff's conjecture]
\label{con:eckhoff}The jump locations of a piecewises-smooth $C^{d}$
function can be reconstructed from its first $\sc$ Fourier coefficients
with asymptotic accuracy $\sc^{-\ord-2}$.
\end{conjecture}
In \cite{batyomAlgFourier} we proposed a reconstruction method (see
\prettyref{alg:half-algo} below) which is based on the original Eckhoff's
procedure.

\begin{algorithm}[H]
\begin{raggedright}
Let $f\in PC\left(\ord,\np\right)$, and assume that $\fun=\sing^{\left(\ord\right)}+\smooth$
where $\sing^{\left(\ord\right)}$ is the piecewise polynomial absorbing
all discontinuities of $\fun$, and $\smooth\in C^{\ord}.$ Assume
in addition the following a-priori bounds:
\par\end{raggedright}
\begin{enumerate}
\item Minimal separation, $\min_{i\neq j}\left|\jp_{i}-\jp_{j}\right|\geq J>0$. 
\item Upper bound on jump magnitudes, $\left|\jc_{l,j}\right|\leq A<\infty$. 
\item Lower bound on the value of the lowest-order jump, $\left|\jc_{0,j}\right|\geq B>0$.
\item Upper bound on the size of the Fourier coefficients of $\smooth$,
$\left|\fc\left(\smooth\right)\right|\leq R\cdot k^{-\ord-2}$. 
\end{enumerate}
\begin{raggedright}
Let us be given the first $3\sc$ Fourier coefficients of $\fun$
for $\sc>\sc\left(\ord,\np,J,A,B,R\right)$ (a quantity which is computable).
The reconstruction is as follows.
\par\end{raggedright}
\begin{enumerate}
\item Obtain first-order approximations to the jump locations $\left\{ \jp_{1},\dots,\jp_{\np}\right\} $
by Prony's method (Eckhoff's method of order 0).
\item Localize each discontinuity $\jp_{j}$ by calculating the first $\sc$
Fourier coefficients of the function $f_{j}=\fun\cdot h_{j}$ where
$h_{j}$ is a $C^{\infty}$ bump function satisfying

\begin{enumerate}
\item $h_{j}\equiv0$ on the complement of $\left[\jp_{j}-J,\jp_{j}+J\right]$;
\item $h_{j}\equiv1$ on $\left[\jp_{j}-\frac{J}{3},\jp_{j}+\frac{J}{3}\right]$.
\end{enumerate}
\item Fix the reconstruction order $\ord_{1}\leq\left\lfloor \frac{\ord}{2}\right\rfloor $.
For each $j=1,2,\dots,\np$, recover the parameters $\left\{ \jp_{j},\jc_{0,j},\dots,\jc_{\ord_{1},j}\right\} $
from the $\ord_{1}+2$ equations
\[
\fc\left(f_{j}\right)=\frac{1}{2\pi}\ee^{-\imath\jp_{j}k}\sum_{l=0}^{\ord_{1}}\frac{\jc_{l,j}}{\left(\imath k\right)^{l+1}},\qquad k=\sc-\ord_{1}-1,\sc-\ord_{1},\dots,\sc
\]
by Eckhoff's method for one jump (in this case we get a single polynomial
equation $\left\{ p_{\sc}^{\ord_{1}}\left(\jp_{j}\right)=0\right\} $
of degree $\ord_{1}$).
\item From the previous steps we obtained approximate values for the parameters
$\left\{ \nn{\jp_{j}}\right\} $ and $\left\{ \nn{\jc}_{l,j}\right\} $.
The final approximation is taken to be
\[
\begin{split}\widetilde{\fun} & =\widetilde{\smooth}+\widetilde{\sing}=\sum_{\left|k\right|\leq\sc}\left\{ \fc(\fun)-\frac{1}{2\pi}\sum_{j=1}^{\np}\ee^{-\imath\nn{\jp_{j}}k}\sum_{l=0}^{d_{1}}\frac{\widetilde{\jc}_{l,j}}{(\imath k)^{l+1}}\right\} \ee^{\imath kx}+\sum_{j=1}^{\np}\sum_{l=0}^{d_{1}}\nn{\jc}_{l,j}V_{l}(x;\nn{\jp_{j}}).\end{split}
\]

\end{enumerate}
\caption{Half-order algorithm, \cite{batyomAlgFourier}.}
\label{alg:half-algo}
\end{algorithm}

We have also shown that this method achieves the following accuracy.
\begin{theorem}[\cite{batyomAlgFourier}]
Let $\fun\in PC\left(\ord,\np\right)$ and let $\nn{\fun}$ be the
approximation of order $\ord_{1}\leq\left\lfloor \frac{\ord}{2}\right\rfloor $
computed by \prettyref{alg:half-algo}. Then%
\footnote{The last (pointwise) bound holds on ``jump-free'' regions.%
}
\begin{equation}
\begin{split}\left|\nn{\jp_{j}}-\jp_{j}\right| & \leq C_{1}\left(\ord,\ord_{1},\np,J,A,B,R\right)\cdot\sc^{-\ord_{1}-2};\\
\left|\nn{\jc}_{l,j}-\jc_{l,j}\right| & \leq C_{2}\left(\ord,\ord_{1},\np,J,A,B,R\right)\cdot\sc^{l-\ord_{1}-1},\; l=0,1,\dots,\ord_{1};\\
\left|\nn{\fun}\left(x\right)-\fun\left(x\right)\right| & \leq C_{3}\left(\ord,\ord_{1},\np,J,A,B,R\right)\cdot\sc^{-\ord_{1}-1}.
\end{split}
\label{eq:half-estimates}
\end{equation}

\end{theorem}
The non-trivial part of the proof of this result was to analyze in
detail the polynomial equation $p\left(\jp_{j}\right)=0$ in step
3 of \prettyref{alg:half-algo}. It turned out that additional orders
of smoothness (namely, between $\ord_{1}$ and $\ord$) produced an
error term which, when substituted into the polynomial $p$, resulted
in unexpected cancellations due to which the root $\jp_{j}$ was perturbed
only by $O\left(\sc^{-\ord_{1}}\right)$. This phenomenon was first
noticed by Eckhoff himself in \cite{eckhoff1995arf} for $\ord=1$,
but at the time its full significance was not realized.

An important property of \prettyref{alg:half-algo} is that its final
asymptotic convergence order essentially depends on the accuracy of
step 3. It is sufficient therefore to replace this step with another
method which achieves full accuracy (i.e. $\sim\sc^{-\ord-2}$) in
order to obtain the overall reconstruction with full accuracy. It
turns out that taking instead of consecutive Fourier samples
\[
k=\sc-\ord-1,\sc-\ord,\dots,\sc
\]
the ``decimated'' section
\begin{eqnarray*}
k & = & \scc,2\scc,\dots,\left(\ord+2\right)\scc;\qquad\scc\isdef\left\lfloor \frac{\sc}{\left(\ord+2\right)}\right\rfloor 
\end{eqnarray*}
provides this accuracy.

For the full details, see \cite{batFullFourier}. Here let us outline
our method of proof.

Denote the single jump point by $\jp\in\left[-\pi,\pi\right]$, and
let $\w=\ee^{-\imath\jp}$.The purpose is to recover the jump point
$\w$ and the jump magnitudes $\left\{ \jc_{0},\dots,\jc_{\ord}\right\} $
from the noisy measurements
\begin{equation}
\widetilde{\fc}\left(\fun\right)=\underbrace{\frac{\w^{k}}{2\pi}\sum_{j=0}^{\ord}\frac{\jc_{j}}{\left(\imath k\right)^{j+1}}}_{\isdef\fc}+\epsilon_{k},\quad k=\scc,2\scc,\dots,\left(\ord+2\right)\scc,\quad\left|\epsilon_{k}\right|\leq R\cdot k^{-\ord-2}.\label{eq:noisy-cf}
\end{equation}
Again, we multiply \eqref{eq:noisy-cf} by $\left(2\pi\right)\left(\imath k\right)^{\ord+1}$.
Denote $\jcc_{j}=\imath^{\ord+1-j}\jc_{\ord-j}$. We get
\begin{equation}
\apprmeas\isdef2\pi\left(\imath k\right)^{\ord+1}\widetilde{\fc}=\underbrace{\w^{k}\sum_{j=0}^{\ord}\jcc_{j}k^{j}}_{\isdef\meas}+\er,\qquad k=\scc,2\scc,\dots,\left(\ord+2\right)\scc,\quad\left|\er\right|\leq R\cdot k^{-1}.\label{eq:noisy-measurements}
\end{equation}

\begin{definition}
Let
\[
p_{\scc}^{\ord}\left(u\right)\isdef\sum_{j=0}^{\ord+1}\left(-1\right)^{j}{\ord+1 \choose j}\meas[\left(j+1\right)\scc]u^{\ord+1-j}.
\]
\end{definition}
\begin{proposition}
The point $u=\w^{\scc}$ is a root of $p_{\scc}^{\ord}\left(u\right)$.
\end{proposition}
\begin{minipage}[t]{1\columnwidth}%
\end{minipage}
\begin{proposition}
The vector of exact magnitudes $\left\{ \jcc_{j}\right\} $ satisfies
\begin{equation}
\begin{bmatrix}\meas[\scc]\w^{-\scc}\\
\meas[2\scc]\w^{-2\scc}\\
\vdots\\
\meas[\left(\ord+1\right)\scc]\w^{-\left(\ord+1\right)\scc}
\end{bmatrix}=\begin{bmatrix}1 & \scc & \scc^{2} & \dots & \scc^{\ord}\\
1 & 2\scc & \left(2\scc\right)^{2} & \dots & \left(2\scc\right)^{\ord}\\
\vdots & \vdots & \vdots & \vdots & \vdots\\
1 & \left(\ord+1\right)\scc & \left(\left(\ord+1\right)\scc\right)^{2} & \dots & \left(\left(\ord+1\right)\scc\right)^{\ord}
\end{bmatrix}\begin{bmatrix}\jcc_{0}\\
\jcc_{1}\\
\vdots\\
\jcc_{\ord}
\end{bmatrix}.\label{eq:the-linear-system}
\end{equation}

\end{proposition}
The procedure for recovery of the $\left\{ \jcc_{0},\dots,\jcc_{\ord},\w\right\} $
is presented in \prettyref{alg:new-reconstruction-single-jump} below,
while the method for full recovery of the function is outlined in
\prettyref{alg:full-algo} below.

\begin{algorithm}[H]
\begin{raggedright}
\begin{minipage}[t]{1\columnwidth}%
Let there be given the first $\scc\gg1$ Fourier coefficients of the
function $\fun_{j}$, and assume that the jump position $\jp$ is
already known with accuracy $o\left(\scc^{-1}\right)$.
\begin{enumerate}
\item Construct the polynomial 
\begin{eqnarray*}
q_{\scc}^{\ord}\left(u\right) & = & \sum_{j=0}^{\ord+1}\left(-1\right)^{j}{\ord+1 \choose j}\apprmeas[\left(j+1\right)\scc]u^{\ord+1-j}.
\end{eqnarray*}
 from the given noisy measurements $\apprmeas[\scc],\apprmeas[2\scc],\dots,\apprmeas[\left(\ord+2\right)\scc]$
\eqref{eq:noisy-measurements}.
\item Find the root $\widetilde{z}$ which is closest to the unit circle
(in fact any root will do).
\item Take $\widetilde{\w}=\sqrt[\scc]{\widetilde{z}}$. Note that in general
there are $\scc$ possible values on the unit circle (see \prettyref{rem:multiple-solutions})
, but since we already know the approximate location of $\w$ the
correct value can be chosen consistently.
\item Recover $\widetilde{\jp}=-\arg\widetilde{\omega}$.
\item To recover the magnitudes, solve the linear system \eqref{eq:the-linear-system}.\end{enumerate}
\begin{remark}
\label{rem:multiple-solutions}To see that there are $\scc$ possible
solutions, notice that one recovers $\ee^{\imath\jp\scc}=\ee^{\imath t}$,
which is satisfied by any $\jp$ of the form $\jp=\frac{t}{\scc}+\frac{2\pi}{\scc}n,\; n\in\ZZ$
and not just $\jp=\frac{t}{\scc}$.\end{remark}
\end{minipage}
\par\end{raggedright}

\caption{Recovery of single jump parameters}
\label{alg:new-reconstruction-single-jump}
\end{algorithm}

\begin{algorithm}[H]
\begin{raggedright}
Let $f\in PC\left(\ord,\np\right)$, and assume that $\fun=\sing^{\left(\ord\right)}+\smooth$
where $\sing^{\left(\ord\right)}$ is the piecewise polynomial absorbing
all discontinuities of $\fun$, and $\smooth\in C^{\ord}.$ Assume
the a-priori bounds as in \prettyref{alg:half-algo}.
\par\end{raggedright}
\begin{enumerate}
\item Using \prettyref{alg:half-algo}, obtain approximate values of the
jumps up to accuracy $O\left(\scc^{-\left\lfloor \frac{\ord}{2}\right\rfloor -2}\right)=o\left(\scc^{-1}\right)$,
and the Fourier coefficients of the functions $\fun_{j}$.
\item Use \prettyref{alg:new-reconstruction-single-jump} to further improve
the accuracy of reconstruction.
\end{enumerate}
\caption{Full accuracy Fourier approximation}
\label{alg:full-algo}
\end{algorithm}

We have shown that indeed full accuracy is acheived.
\begin{theorem}[\cite{batFullFourier}]
\label{thm:jump-accuracy-new}\prettyref{alg:new-reconstruction-single-jump}
recovers the parameters of a single jump from the given noisy measurements
\eqref{eq:noisy-measurements} with the following accuracy:
\begin{eqnarray*}
\left|\nn{\jp}-\jp\right| & \leq & C_{4}\frac{R}{B}\scc^{-\ord-2},\\
\left|\nn{\jcc}_{j}-\jcc_{j}\right| & \leq & C_{5}R\left(1+\frac{A}{B}\right)\scc^{-j-1},\qquad j=0,1,\dots,\ord.
\end{eqnarray*}

\end{theorem}
The main idea of the proof is to analyze the perturbation of the polynomial
$p_{\scc}^{\ord}\left(u\right)$ by $q_{\scc}^{\ord}$ using Rouche's
theorem.

After making the substitution $N=\left\lfloor \frac{\sc}{\left(\ord+2\right)}\right\rfloor $,
we obtain as an immediate consequence of \prettyref{thm:jump-accuracy-new}
the resolution of \prettyref{con:eckhoff}.
\begin{theorem}
\label{thm:full-accuracy-final}Let $\fun\in PC\left(\ord,\np\right)$
and let $\nn{\fun}$ be the approximation of order $\ord$ computed
by \prettyref{alg:full-algo}. Then
\begin{equation}
\begin{split}\left|\nn{\jp_{j}}-\jp_{j}\right| & \leq C_{6}\left(\ord,\np,J,A,B,R\right)\cdot\sc^{-\ord-2};\\
\left|\nn{\jc}_{l,j}-\jc_{l,j}\right| & \leq C_{7}\left(\ord,\np,J,A,B,R\right)\cdot\sc^{l-\ord-1},\qquad l=0,1,\dots,\ord;\\
\left|\nn{\fun}\left(x\right)-\fun\left(x\right)\right| & \leq C_{8}\left(\ord,\np,J,A,B,R\right)\cdot\sc^{-\ord-1}.
\end{split}
\label{eq:full-estimates}
\end{equation}

\end{theorem}
Note that the system \eqref{eq:noisy-cf} is a certain variant of
the confluent Prony system \eqref{eq:confluent equation_prony_system}
for just one node. Therefore, \prettyref{alg:new-reconstruction-single-jump}
can be regarded as a concrete solution method for this particular
case.

\bibliographystyle{plain}
\bibliography{../../../bibliography/all-bib}

\end{document}